%% file: main.tex
\documentclass[10pt]{article}
\usepackage{graphicx}
\usepackage{cite}
\usepackage{natbib}
\usepackage{todonotes}
\usepackage{algpseudocode}
\usepackage{algorithm2e}
\RestyleAlgo{ruled}
\usepackage{amsmath}
\usepackage{amssymb}
\usepackage{amsthm}
\usepackage{gensymb}
\usepackage[right]{lineno}
\usepackage{authblk}

\usepackage{multirow}
\usepackage{mathrsfs}

\usepackage{geometry}
\RequirePackage[OT1]{fontenc}
\RequirePackage[colorlinks,citecolor=blue,urlcolor=blue]{hyperref}

\newcommand{\abs}[1]{\left \vert #1 \right \vert}

\newcommand{\cT}[1]{\mathcal{#1}}
\newcommand{\bT}[1]{\mathbb{#1}}

\newcommand{\RM}[1]{\mathrm{#1}}

\newcommand{\norm}[1]{\left\lVert#1\right\rVert}

\newcommand{\BT}[1]{\textbf{#1}}

\numberwithin{equation}{section}
\theoremstyle{plain}
\newtheorem{theorem}{Theorem}[section]
\newtheorem{lemma}{Lemma}[section]
 \newtheorem{corollary}{Corollary}[section]

\theoremstyle{definition}
 \newtheorem{definition}{Definition}[section]
  \newtheorem{notation}{Notation}[section]
    
    \newtheorem{remark}{Remark}[section]

\usepackage{siunitx}
\usepackage{lineno} 
\usepackage{setspace} 
\usepackage{threeparttable}
\usepackage{booktabs}
\usepackage{hyperref}
\usepackage{color}
\usepackage{soul}

\usepackage[normalem]{ulem}
\usepackage{textcomp}

\newcommand{\txt}[1]{\iffalse #1 \fi}
\newcommand{\red}[1]{\textcolor{red}{#1}}

\renewcommand{\txt}[1]{\textcolor{red}{#1}}

\begin{document}

\numberwithin{equation}{section}

\author[1,2, *]{N. Alexia Raharinirina}
\author[1,3]{Konstantin Fackeldey}
\author[1]{Marcus Weber}
\affil[1]{\small Zuse Institute Berlin, Takustra{\ss}e 7, 14195 Berlin, Germany}
\affil[2]{\small Free University of Berlin, Kaiserswerther Str. 16-18, 14195 Berlin}
\affil[3]{\small Technical University Berlin, Institute for Mathematics, Stra{\ss}e des 17. Juni 136, 10623 Berlin, Germany}
\affil[*]{corresponding author: n.raharinirina@fu-berlin.de, missralexia1@gmail.com}

\title{Qualitative Euclidean embedding of Disjoint Sets of Points}

\hypersetup{pdftitle={Euclidean embedding of separated sets of points}}

\maketitle

\vspace*{-1cm}
\begin{abstract}               
\input{abstract}

\textbf{keywords}: Euclidean embedding, Ger\v{s}gorin circle theorem, Gerschgorin circle theorem, Schoenberg criterion, disjoint sets, Euclidean distance matrix
\end{abstract}

\section*{Declarations}
\input{Declarations}

\section{Prerequisites}
\input{prerequisite}

\section{Main result}
\input{finding}

\section{Intermediate results}
\label{sec:details}
\input{intermediate}

\bibliographystyle{cell}
\bibliography{main}
\end{document}

%% file: abstract.tex
We consider two disjoint sets of points. If at least one of the sets can be embedded into an Euclidean space, then we provide sufficient conditions for the two sets to be jointly embedded in one Euclidean space. In this joint Euclidean embedding, the distances between the points are generated by a specific relation-preserving function. Consequently, the mutual distances between two points of the same set are specific qualitative transformations of their mutual distances in their original space; the pairwise distances between the points of different sets can be constructed from an arbitrary proximity function. 

%% file: Declarations.tex
\paragraph{Funding:} 
Funded by the Deutsche Forschungsgemeinschaft (DFG, German Research Foundation) under Germany's Excellence Strategy – The Berlin Mathematics Research Center MATH+ (EXC-2046/1, project ID: 390685689), project EF5-4 ``The Evolution of Ancient Egyptian – Quantitative and Non-Quantitative Mathematical Linguistics'' and the Einstein Center Chronoi, project ``Diglossie im ptolem{\"a}ischen {\"A}gypten'' (FU Berlin).

\paragraph{Code availability:}
To demonstrate our results, python codes and jupyter notebooks are provided at:\\ \url{https://github.com/AlexiaNomena/PSD_cosine_law_matrix}






%% file: prerequisite.tex
In this section, we provide the most important concepts and theorems that are necessary for understanding our proofs. The general notation scheme is shown in Table~\ref{tab:Nom}.

\begin{table}[ht]
    \centering
    \caption{\textbf{General notation scheme.}}
    \small
    \begin{tabular}{@{}p{0.4\linewidth}p{0.4\linewidth}p{0.2\linewidth}@{}}
    \toprule
        Objects & Notations & Examples\\
        \midrule
             Scalar / integer & italic upper case letter or lower case letters  & $i$, $m$, $n$, $I$, $M$, $N, \dots$\\
            Set & curly upper case letters  & $\cT{X},\,\cT{Y}, \cT{V}, \dots$\\
            Distance measure in sets & $d$ with the set in subscript  & $d_{\cT{X}},d_{\cT{Y}}, d_{\cT{V}}, \dots$\\
            Point of a set $\RM{v}_s$ & straight lower case & $\RM{x}_m$, $\RM{y}_n$,  $\RM{v}_s, \dots$\\
            Euclidean embedding of a point $\RM{v}_s$ & same notation as the point & $\RM{x}_m$, $\RM{y}_n$,  $\RM{v}_s, \dots$\\
            Matrix  & bold upper case letter  & $\BT{A}, \, \BT{B}$, $\BT{C},\dots$ \\
            Dot product & a point-sign between vectors/matrices & $\RM{x}\cdot \RM{y}$, $\RM{x}\cdot \BT{A}$, $\BT{A} \cdot \BT{B}$\\
            Row vector pointing from an origin to a point & same notation as the point & $\RM{x}_m$, $\RM{y}_n$,  $\RM{v}_s, \dots$\\
            A place holder for variable items & a sign in the position of the variable & $\bullet$ \\
            Dependency to variables (case 1) & anything with parenthesis & $f(\bullet,\, \bullet)$, $\BT{A}(\bullet), \dots$\\
            Dependency to variables (case 2) & subscripts and superscripts  & $i_{\bullet}$, $\BT{A}_{\bullet}$, $\cT{P}^{\bullet}_{\bullet}, \dots$\\
            An element of a matrix $\BT{A}$ & bold lower case with subscripts     & $\BT{a}_{m, n}$ \\
            Euclidean norm & norm notation with the space as subscript & $\norm{\bullet}_{\bT{R}^{\bullet}}$\\
    \bottomrule
    \end{tabular}
    \label{tab:Nom}
\end{table}

\begin{definition}{(Distance set)}
\label{def:dist}

Let $\cT{V}$ be an arbitrary a set of $S\in \bT{N}\setminus \{0, 1\}$ points
\begin{equation}
\cT{V}:=\{\RM{v}_s: s \in \{1,\ldots,\,S\}\},
\label{eq:suppV}
\end{equation}
and let $d_{\cT{V}}:\cT{V}\times \cT{V} \rightarrow \, \bT{R}^+ $ be a distance function associated with $\cT{V}$.

We consider the set of numbers
\begin{equation}
\label{eq:suppK}
\cT{D}_{\cT{V}} := \Big\{d_{\cT{V}}(\RM{v}_l, \RM{v}_s) = d_{\cT{V}}(\RM{v}_s, \RM{v}_l): \RM{v}_l, \RM{w}_s \in \cT{V} \Big\},
\end{equation}
containing the pairwise distance between the points of $\cT{V}$. We call $\cT{D}_{\cT{V}}$ the ``distance set'' associated with $\cT{V}$.
\end{definition}

\begin{definition}{(Proximity function, proximity set, and symmetric proximity level)}
\label{def:prox}

A ``proximity function'' between two sets of points $\cT{V}$ and $\cT{W}$ is a non-negative function that maps any pairs of points to some measure of closeness between them. $f$ is characterized as follows
\begin{equation}
\begin{split}
    f:\quad \cT{V} \times \cT{W} &  \rightarrow \, \bT{R}^+,\\
      \, (\RM{v},\,\RM{w}) &  \mapsto f(\RM{v},\, \RM{w}) \geq 0.
\end{split}
\end{equation}
We then define the associated the ``proximity set'' as
\begin{equation}
\label{eq:suppF}
\cT{P}^f_{\cT{V}\times \cT{W}} := \Big\{p^f_{\cT{V}\times \cT{W}}(\RM{v},\,\RM{w}) = p^f_{\cT{V} \times \cT{W}}(\RM{w},\,\RM{v}) = f(\RM{v},\,\RM{w}) : \RM{v} \in \cT{V},\,\RM{w} \in \cT{W} \Big\}.
\end{equation}
we refer to $p^f_{\cT{V}\times \cT{W}}$ as the ``symmetric proximity level'' corresponding to $f$. If $\cT{W} = \cT{V}$ we simply denote the proximity set $\cT{P}^f_{\cT{V}}$ and the symmetric proximity level $p^f_{\cT{V}}$ in Eq.~\ref{eq:suppF} (and vice-versa if $\cT{V} = \cT{W}$).
\end{definition}

\begin{definition}{(Euclidean embedding of a set of points)}
\label{def:embeucl}

Let $\cT{V}$ be an arbitrary set of $S\in \bT{N}\setminus \{0, 1\}$ points with associated distance set $\cT{D_{\cT{V}}}$ as defined in Def.~\ref{def:dist} (resp. proximity set $\cT{P}^f_{\cT{V}}$).

We say that $\cT{V}$ can be embedded into the Euclidean space $\bT{R}^Q$ with respect to the proximity set $\cT{P}^{d_{\cT{V}}}_{\cT{V}}$ (resp. $\cT{P}^f_{\cT{V}}$) it satisfies the following 
\begin{itemize}
    \item[(i)] there exists a bijection $\phi$ between $\cT{V}$ and a subspace of $\bT{R}^Q$.
    \item[(ii)] for any pair of points $(\RM{v}_l,\RM{v}_s) \in \cT{V}\times\cT{V}$,
    \begin{equation*}
       \norm{\phi(\RM{v}_l) - \phi(\RM{v}_s)}_{\bT{R}^Q} = d_{\cT{V}}(\RM{v}_l,\RM{v}_s) \quad \text{(resp. } \quad  \norm{\phi(\RM{v}_l) - \phi(\RM{v}_s)}_{\bT{R}^Q} = p^f_{\cT{V}}(\RM{v}_l,\RM{v}_s)\text{)} 
    \end{equation*}
\end{itemize}
\end{definition}

\begin{definition}{(Set of mutual Euclidean distances)}
\label{def:disteucl}

Let $\cT{V}$ be an arbitrary set of $S\in \bT{N}\setminus \{0, 1\}$ points with associated distance set $\cT{D_{\cT{V}}}$ (resp. proximity set $\cT{P}^f_{\cT{V}}$. 

The set of numbers $\cT{D_{\cT{V}}}$ (resp. $\cT{P}^f_{\cT{V}}$) is a ``set of mutual Euclidean distances'', between the points of $\cT{V}$, if $\cT{V}$ with respect to the proximity set  $\cT{P}^{d_{\cT{V}}}_{\cT{V}}$ (resp. $\cT{P}^f_{\cT{V}}$) can be embedded into an Euclidean space.
\end{definition}

\begin{definition}{(Cosine law matrix)}
\label{def:M}

Let $\cT{V}$ be an arbitrary set of $S\in \bT{N}\setminus \{0, 1\}$ and let $$h:\cT{V}\times \cT{V} \rightarrow \bT{R}^+,$$ a symmetric proximity function or a symmetric proximity level. 

Then, for an arbitrary index $a\in \{1, \ldots, S\}$ let $\BT{M}^{\:\RM{v}_a}(\cT{V}, h)$ be the real symmetric $S\times S$ matrix of components
\begin{equation}
\label{eq:mat}
    \BT{m}_{l, s}^{\:\RM{v}_a}(\cT{V}, h) := \BT{m}_{s, l}^{\:\RM{v}_a}(\cT{V}, h) := \frac{h(\RM{v}_l,\,\RM{v}_a)^2 + h(\RM{v}_a,\,\RM{v}_s)^2 - h(\RM{v}_l,\,\RM{v}_s)^2}{2},
\end{equation}
We call  $\BT{M}^{\:\RM{v}_a}(\cT{V}, h)$ ``cosine law matrix'' associated with $(\cT{V}, h)$ and reference $\RM{v}_a$ because it is based on the cosine law in Euclidean space where $\RM{v}_a$ is assumed to be used as a reference point.
\end{definition}

\begin{definition}{(Positive semi-definite matrix)}

Let $\BT{A}$ be an arbitrary $I\times I$ symmetric real matrix ($I\in \bT{N}\setminus \{0, 1\}$). $\BT{A}$ is positive semi-definite if it satisfies one of the following equivalent statements
\begin{itemize}
    \item[(i)] $\RM{x}\cdot \BT{A}\cdot \RM{x}^{\intercal}$ where $\RM{x}$ is a real row vector of $\bT{R}^I$,
    \item[(ii)] all the eigenvalues of $\BT{A}$ are real and non-negative.
\end{itemize}
\end{definition}

\begin{lemma} 
\label{lm:sumPSD}
Let $\BT{A}_p$ ($p\in \{1,\ldots,\,P\}, P\in \bT{N}\setminus\{0\}$) be positive semi-definite matrices. The matrix
$$\BT{A} = \sum\limits_{p = 1}^{P} \BT{A}_p,$$ is positive semi-definite.
\end{lemma}

\begin{theorem}{(Young-Householder or Schoenberg criterion)}
\label{thm:YH}

Let $\cT{V}$ be an arbitrary set of $S\in \bT{N}\setminus \{0, 1\}$ points with an associated distance set $\cT{D_{\cT{V}}}$ associated with the distance function $d_{\cT{V}}$ (resp. proximity set $\cT{P}^f_{\cT{V}}$ associated with symmetric proximity level $p^f_{\cT{V}}$). 

The set $\cT{D_{\cT{V}}}$ (resp. $\cT{P}^f_{\cT{V}}$) is a set of mutual Euclidean distances, between the points of $\cT{V}$, if and only if the cosine law matrix $\BT{M}^{\:\RM{v}_a}(\cT{V},\, d_{\cT{V}})$ (resp. $\BT{M}^{\:\RM{v}_a}(\cT{V},\, p^f_{\cT{V}})$) is positive semi-definite. This result does not depend on the choice of $\RM{v}_a$.

\end{theorem}
\begin{proof}
See \citep{YH1938} or \citep{Schoenberg1935}
\end{proof}

\begin{corollary}
\label{col:YH}
Let $\cT{V}$ be an arbitrary set of $S\in \bT{N}\setminus \{0, 1\}$ points with an associated distance set $\cT{D_{\cT{V}}}$ (resp. proximity set $\cT{P}^f_{\cT{V}}$. 

The set $\cT{V}$ can be embedded into the Euclidean space $\bT{R}^Q$ with respect to $\cT{D_{\cT{V}}}$ (resp. $\cT{P}^f_{\cT{V}}$), if and only if the associated cosine law matrix can be decomposed as 
$$\BT{M}^{\:\RM{v}_a}(\cT{V},\,\bullet) = \BT{E}\cdot \BT{E}^{\,\intercal},$$
where $\BT{E}$ is a $S \times S$ matrix. The $s^{th}$ row of $\BT{E}$ is the Euclidean embedding of the point $\RM{v}_s$ for all $s \in \{1, \ldots, S\}$.
\end{corollary}
\begin{proof}
By definition~\ref{def:embeucl} $\cT{V}$ can be embedded into the Euclidean space if the distance set $D_{\cT{V}}$ (resp.  proximity set $\cT{P}^f_{\cT{V}}$) is a set of mutual Euclidean distance. This is true if an only if and only if cosine law matrix $\BT{M}^{\:\RM{v}_a}(\cT{V}, \,\bullet)$ is positive semi-definite (Theorem~\ref{thm:YH}). We have two steps for proving this. 

``$\Rightarrow$'': Let us assume that $\BT{M}^{\:\RM{v}_a}(\cT{V}, \,\bullet)$ is positive semi-definite. Since $\BT{M}^{\:\RM{v}_a}(\cT{V}, \,\bullet)$ is real and symmetric, there exists an eigendecomposition such that $$\BT{M}^{\:\RM{v}_a}(\cT{V}, \bullet) = \BT{U}\cdot \BT{D}_{\lambda} \cdot \BT{U}^{\,\intercal},$$
where $\BT{U}$ is a unitary $S\times S$ matrix and $\BT{D}_{\lambda}$ is the diagonal matrix of eigenvalues of $\BT{M}^{\:\RM{v}_a}(\cT{V}, \bullet)$ compatible with $\BT{U}$. Since $\BT{M}^{\:\RM{v}_a}(\cT{V}, \,\bullet)$ is positive semi-definite all the eigenvalues are non-negative and thus we can choose $$\BT{E} = \BT{U}\cdot \BT{D}_{\lambda}^{1/2},$$ where $\BT{D}_{\lambda}^{1/2}$ is the diagonal matrix of square roots of eigenvalues $\BT{M}^{\:\RM{v}_a}(\cT{V}, \bullet)$, to obtain
$$\BT{M}^{\:\RM{v}_a}(\cT{V},\,\bullet) = \BT{E}\cdot \BT{E}^{\,\intercal}.$$

``$\Leftarrow$'': Let $\BT{E}$ be a real $S\times S$ matrix such that 
$$\BT{M}^{\:\RM{v}_a}(\cT{V},\,\bullet) = \BT{E}\cdot \BT{E}^{\,\intercal}.$$
And let $\RM{x}$ be a real row vector of $\bT{R}^S$
\begin{align*}
    \RM{x}\cdot\BT{M}^{\:\RM{v}_a}(\cT{V},\,\bullet)\cdot\RM{x}^{\intercal} & =\, \RM{x}\cdot\BT{E}\cdot \BT{E}^{\,\intercal}\cdot\RM{x}^{\intercal}\\
    &=\, \left(\RM{x}\cdot\BT{E}\right)\cdot \left(\RM{x}\cdot\BT{E}\right)^{\,\intercal} \geq 0
\end{align*}
by definition of the dot product. Thus, $\BT{M}^{\:\RM{v}_a}(\cT{V},\,\bullet)$ is positive semi-definite.
\end{proof}

\begin{theorem}{(Ger\v{s}gorin circle theorem)}
\label{thm:G}
Let $\BT{A}$ be an arbitrary $I\times I$ symmetric complex matrix ($I\in \bT{N}\setminus \{0, 1\}$) and let $R_i$ be the sum of the modulus of the non-diagonal entries in the $i$-th row of $\BT{A}$, that is
\begin{equation*}
 R_i = \sum\limits_{j\neq i} \abs{\BT{a}_{i, j}}.   
\end{equation*}
The closed disc centered at $\BT{a}_{i, i}$ with radius $R_i$, denoted $D(\BT{a}_{i, i},\, R_i) \subsetneq \bT{C}$, is called a Ger\v{s}gorin disc. 

For every eigenvalues $\lambda$ of $\BT{A}$, there exists an index $i_{\lambda}\in \{1, \ldots, I\}$ such that $\lambda \in D(\BT{a}_{i_{\lambda}, i_{\lambda}},\, R_{i_{\lambda}})$. In other words
$$\forall \text{ eigenvalue } \lambda \text{ of } \BT{A}, \exists\: i_{\lambda} \in \{1, \ldots, I\}, \text{ such that } \abs{\lambda - \BT{a}_{i_{\lambda}, i_{\lambda}}} \,\leq\, R_{i_{\lambda}}.$$
\end{theorem}
\begin{proof}
See \citep{Gerschgorin1931} or \citep{Horn1985}.
\end{proof}

\begin{notation}{(Two arbitrary sets of points)}
\label{not:2sets}

Let us consider two arbitrary sets of $M$ and $N$ points ($M, N \in \bT{N}\setminus \{0\}$) denoted
$$\cT{X} := \left\{\RM{x}_m : m \in \{1,\ldots,\,M\}\right\},$$
and 
$$\cT{Y} := \left\{\RM{y}_n : n \in \{1,\ldots,\,N\}\right\}.$$
We denote distance sets $\cT{D}_{\cT{X}}$ and $\cT{D}_{\cT{Y}}$, associated to $\cT{X}$ and $\cT{Y}$, respectively. We also consider some proximity function  $$u_{\cT{X}}: \cT{X}\times \{\RM{o}_{\cT{X}}\} \rightarrow \bT{R}^+ \quad \text{and} \quad u_{\cT{Y}}: \cT{Y}\times \{\RM{o}_{\cT{Y}}\} \rightarrow \bT{R}^+,$$
where $\RM{o}_{\cT{X}}$ and $\RM{o}_{\cT{Y}}$ are points of origin such that the proximity function $u_{\bullet}$ define the distance to origin for each corresponding sets. We denote the associated proximity sets to origin as $\cT{P}^{u_\cT{X}}_{\cT{X}\times \{\RM{o}_{\cT{X}}\}}$ and $\cT{P}^{u_\cT{Y}}_{\cT{Y}\times \{\RM{o}_{\cT{Y}}\}}$. 

All set of numbers $\cT{D}_{\cT{X}}$, $\cT{D}_{\cT{Y}}$, $\cT{P}^{u_\cT{X}}_{\cT{X}\times \{\RM{o}_{\cT{X}}\}}$, and $\cT{P}^{u_\cT{Y}}_{\cT{Y}\times \{\RM{o}_{\cT{Y}}\}}$, must contain at least one non-zero value.

And lastly, we consider 
\begin{itemize}
    \item [(i)] an arbitrary proximity function $f$ measuring some notion of closeness between the points of different sets $\cT{X}$ and $\cT{Y}$, such that $f
    $ is not $0$ everywhere.
    \item [(ii)] theoretical origin point $\RM{o}$ for defining the distance to the origin of the points of $\cT{X}\cup \cT{Y}.$
\end{itemize} 
\end{notation}
\begin{notation}
\label{not:vz}
Given $\cT{X}$ and $\cT{Y}$ as characterized in Notation~\ref{not:2sets}, we denote
$$\cT{W}:=\cT{X}\cup \cT{Y} \cup \{\RM{o}\}.$$
With the following convention for the points $\RM{w}_s$ of $\cT{W}$, where $s \in \{1, \ldots,\, S = M+N+1\}$
\begin{equation}
    \RM{w}_s = 
    \begin{cases}
    \RM{x}_{s} & \quad \text { if } s=1, \ldots, M,\\
     \RM{o} & \quad \text { if } s=M+1,\\
     \RM{y}_{s-(M+1)} & \quad \text { if } s=M+2, \ldots, M+N+1,
    \end{cases}
\end{equation}

We also denote
$$\cT{V}^{\,\RM{z}} := \cT{W} \cup \{\RM{z}\}, $$
where $\RM{z} \notin \cT{X}$ and $\RM{z} \notin \cT{Y}$ is an arbitrary theoretical point. 

We use the following convention for the points $\RM{v}_q$ of $\cT{V}^{\,\RM{z}}$ where $q \in \{1, \ldots,\, Q = M+N+2\}$
\begin{equation}
    \RM{v}_q = 
    \begin{cases}
    \RM{z} & \quad \text { if } q=1,\\
    \RM{x}_{q-1} & \quad \text { if } q=2, \ldots, M+1,\\
     \RM{o} & \quad \text { if } q=M+2,\\
     \RM{y}_{q-(M+2)} & \quad \text { if } q=M+3, \ldots, M+N+2,
    \end{cases}
\end{equation}
\end{notation}

%% file: finding.tex
In this section, we provide the main result of our paper. 
\begin{theorem}{(Sufficient conditions for a Qualitative Euclidean embedding)}
\label{thm:main}

Let us consider the set of $M$ and $N$ points, $\cT{X}$ and $\cT{Y}$, with corresponding distance sets $\cT{D}_{\cT{X}}$ and $\cT{D}_{\cT{Y}}$, with a defined proximity sets to origin $\cT{P}^{u_{\cT{X}}}_{\cT{X}\times \{\RM{o}_{\cT{X}}\}}$ and $\cT{P}^{u_{\cT{Y}}}_{\cT{Y}\times \{\RM{o}_{\cT{Y}}\}}$, an arbitrary proximity function $f$, and a theoretical origin point $\RM{o}$, as characterized in Notation~\ref{not:2sets}. With the notation $\cT{W} = \cT{X}\cup \cT{Y} \cup \{\RM{o}\}$, with $\cT{W}$ characterized in Notation~\ref{not:vz}.

Then, there exists $\varepsilon > 0$ parameterizing a proximity function $f^{\varepsilon}$ given by,
\begin{equation}
\label{eq:f_varepsilon}
\begin{split}
  f^{\varepsilon}: \cT{W}\times \cT{W} \rightarrow \bT{R}^+
\end{split}
\end{equation}
satisfying that for every $\RM{w}_{s_1},\,\RM{w}_{s_2}\in \cT{W}$ where $s_1, s_2\in \{1, \ldots S = M+N+1\}$,
\begin{equation}
    f^{\,\varepsilon}(\RM{w}_{s_1},\,\RM{w}_{s_2}) =
    \begin{cases}
          0 & \text{ if } \RM{w}_{s_1} = \RM{w}_{s_2}, \\
          
          \left(p^{u_{\cT{A}}}_{\cT{A}\times\{\RM{o}_{\cT{A}}\}}(\RM{w}_{s_1},\,\RM{o}_{\cT{A}})^2 + \varepsilon\right)^{1/2}
          & \text{ if } \RM{w}_{s_1} \in \cT{A},\,\RM{w}_{s_2} = \RM{o},\, \cT{A} = \cT{X} \text{ or } \cT{Y},\\
          
          \left(p^{u_{\cT{A}}}_{\cT{A}\times\{\RM{o}_{\cT{A}}\}}(\RM{w}_{s_2},\,\RM{o}_{\cT{A}})^2 + \varepsilon\right)^{1/2}
          & \text{ if } \RM{w}_{s_1}  = \RM{o},\,\RM{w}_{s_2}\in \cT{A},\, \cT{A} = \cT{X} \text{ or } \cT{Y},\\
          
          \left(d_{\cT{A}}(\RM{w}_{s_1},\,\RM{w}_{s_2})^2 + \varepsilon\right)^{1/2} & \text{ if } \RM{w}_{s_1},\, \RM{w}_{s_2}\in \cT{A},\, \cT{A} = \cT{X} \text{ or } \cT{Y},\\
          
          \left(p^f_{\cT{X}\times \cT{Y}}(\RM{w}_{s_1},\,\RM{w}_{s_2})^2 + \varepsilon \right)^{1/2} 
          & \text{ if } \RM{w}_{s_1} \in \cT{X}, \, \RM{w}_{s_2}\in \cT{Y} \text{ or vice-versa}.
    \end{cases}
\end{equation}
And furthermore, the set of points $\cT{W}$, with respect to the proximity set $\cT{P}^{f^{\varepsilon}}_{\cT{W}}$, can be embedded into the Euclidean space $\bT{R}^Q$ ($Q = M+N+2$).
\end{theorem}
\begin{remark}
Here, we note that Theorem~\ref{thm:main} is also valid without including the theoretical point of origin $\RM{o}$.
\end{remark}

\begin{proof}[proof of Theorem~\ref{thm:main}]
Let $\alpha\geq 0$ and let us consider the symmetric proximity function $f^{\,\alpha}$ given by
$$f^{\,\alpha}:{\cT{W}\times\cT{W}} \rightarrow \bT{R}^+,$$
satisfying that, for every $\RM{w}_{s_1},\,\RM{w}_{s_2}\in \cT{W}$ where $s_1, s_2\in \{1, \ldots S = M+N+1\}$
\begin{equation}
\label{eq:f_alpha}
    f^{\,\alpha}(\RM{w}_{s_1},\,\RM{w}_{s_2}) =
    \begin{cases}
          0 & \text{ if } \RM{w}_{s_1} = \RM{w}_{s_2}, \\
          
          \left(p^{u_{\cT{A}}}_{\cT{A}\times\{\RM{o}_{\cT{A}}\}}(\RM{w}_{s_1},\,\RM{o}_{\cT{A}})^2 + \alpha\right)^{1/2}
          & \text{ if } \RM{w}_{s_1} \in \cT{A},\,\RM{w}_{s_2} = \RM{o},\, \cT{A} = \cT{X} \text{ or } \cT{Y},\\
          
          \left(p^{u_{\cT{A}}}_{\cT{A}\times\{\RM{o}_{\cT{A}}\}}(\RM{w}_{s_2},\,\RM{o}_{\cT{A}})^2 + \alpha\right)^{1/2}
          & \text{ if } \RM{w}_{s_1}  = \RM{o},\,\RM{w}_{s_2}\in \cT{A},\, \cT{A} = \cT{X} \text{ or } \cT{Y},\\
          
          \left(d_{\cT{A}}(\RM{w}_{s_1},\,\RM{w}_{s_2})^2 + \alpha\right)^{1/2} & \text{ if } \RM{w}_{s_1},\, \RM{w}_{s_2}\in \cT{A},\, \cT{A} = \cT{X} \text{ or } \cT{Y},\\
          
          \left(p^f_{\cT{X}\times \cT{Y}}(\RM{w}_{s_1},\,\RM{w}_{s_2})^2 + \alpha \right)^{1/2} 
          & \text{ if } \RM{w}_{s_1} \in \cT{X}, \, \RM{w}_{s_2}\in \cT{Y} \text{ or vice-versa} ,
    \end{cases}
\end{equation}
where $p^f_{\cT{X}\times \cT{Y}}$ is the symmetric proximity level $\cT{P}^f_{\cT{X}\times \cT{Y}}$ associated with $f$, and similarly $p^{u_{\cT{A}}}_{\cT{A}\times\{\RM{o}_{\cT{A}}\}}(\bullet, \bullet)$ is the proximity level associated with $u_{\cT{A}}$ ($\cT{A} = \cT{X} \text{ or } \cT{Y}$). 

Let us also consider a theoretical point $\RM{z}$ and the set of point $\cT{V}^{\,\RM{z}} = \cT{W}\cup \{z\}$ (Notation~\ref{not:vz}). We then provide $\cT{V}^{\,\RM{z}}$ with the following proximity set 
\begin{equation}
\label{eq:proxVz}
\cT{P}^{h^{\,\alpha}}_{\cT{V}^{\,\RM{z}}} := \Big\{p^{h^{\,\alpha}}_{\cT{V}^{\,\RM{z}}}(\RM{v}_{q_1},\,\RM{v}_{q_2}) = p^{h^{\,\alpha}}_{\cT{V}^{\,\RM{z}}}(\RM{v}_{q_1},\,\RM{v}_{q_2}) = {h^{\,\alpha}}(\RM{v}_{q_1},\,\RM{v}_{q_2}): \RM{v}_{q_1},\,\RM{v}_{q_2} \in \cT{V}^{\,\RM{z}}\Big\},
\end{equation}
where $q_1, q_2 \in \{1, \ldots, Q = M+N+2\}$ and $h^{\,\alpha}:\cT{V}^{\,\RM{z}} \times \cT{V}^{\,\RM{z}} \rightarrow \bT{R}^+$ is the symmetric proximity function
\begin{equation}
\label{eq:h_alpha}
    h^{\,\alpha}(\RM{v}_{q_1},\,\RM{v}_{q_2}) =
    \begin{cases}
    0 & \text{ if } \RM{v}_{q_1} = \RM{v}_{q_2}, \\
    
    f^{\,\alpha}(\RM{v}_{q_1},\,\RM{v}_{q_2}) & \RM{v}_{q_1} \neq \RM{v}_{q_2} \text{ and } \RM{v}_{q_1},\,\RM{v}_{q_2}\neq \RM{z},\\
    
    p_{\cT{V} \times \{\RM{z}\}}^{\, g}(\RM{v}_{q_1},\,\RM{z}) & \text{ if } \RM{v}_{q_2} = \RM{z}, \text{ and } \RM{v}_{q_1}\neq \RM{v}_{q_2},\\
    
    p_{\cT{V} \times \{\RM{z}\}}^{\, g}(\RM{v}_{q_2},\,\RM{z}) & \text{ if } \RM{v}_{q_1} = \RM{z} \text{ and } \RM{v}_{q_1}\neq \RM{v}_{q_2}.
    \end{cases}
\end{equation}
where $g:\cT{W} \times \{\RM{z}\} \rightarrow \bT{R}^+$ is an arbitrary proximity function that does not necessarily depend on $\alpha$ and $p_{\cT{V} \times \{\RM{z}\}}^{\, g}$ is the symmetric proximity level corresponding to $g$ (Def.~\ref{def:prox}). 

It is evident that, for any arbitrary point $\RM{z}$, if $\cT{V}^{\,\RM{z}}$ can be embedded in the Euclidean space $\bT{R}^Q$, then $\cT{W}$ can also be embedded in $\bT{R}^Q$. Therefore, we search for a sufficient condition to obtain an Euclidean embedding of $\cT{V}^{\,\RM{z}}$. This Euclidean embedding can be found when the proximity set $\cT{P}^{h^{\,\alpha}}_{\cT{V}^{\,\RM{z}}}$ is a set of mutual Euclidean distances between the points of $\cT{V}^{\,\RM{z}}$.

Firstly, according to Theorem~\ref{thm:YH}, the proximity set $\cT{P}^{h^{\,\alpha}}_{\cT{V}^{\,\RM{z}}}$ is a set of mutual Euclidean distance, between the points of $\cT{V}^{\,\RM{z}}$, if and only if the cosine law matrix $\BT{M}^{\:\RM{z}}(\cT{V}^{\,\RM{z}},\, h^{\,\alpha})$ is positive semi-definite. 

Next, we need to find a sufficient condition for $\BT{M}^{\:\RM{z}}(\cT{V}^{\,\RM{z}},\, h^{\,\alpha})$ to be positive semi-definite and the dimensions Euclidean embedding space is equal to the number of rows/columns of  $\BT{M}^{\:\RM{z}}(\cT{V}^{ \,\RM{z}},\, h^{\,\alpha})$ which is equal to $Q$ (see corollary~\ref{col:YH}). This is tackled in section~\ref{sec:details} (notably in Theorem~\ref{thm:vareps}) and relies on a specific construction of the point $\RM{z} = \RM{z}^*$ and on a specific choice of $\alpha = \varepsilon >0 $. From this procedure, we conclude that there exists a point $\RM{z}^*$ and $\varepsilon>0$  such that $\BT{M}^{\:\RM{z}^*}(\cT{V}^{\,\RM{z}^*},\, h^{\varepsilon})$ is positive semi-definite and thus $\cT{V}^{\,\RM{z}^*}$ can be embedded in the Euclidean space $\bT{R}^Q$, and thus, $\cT{W}$ can also be embedded in $\bT{R}^Q$. Finally, the set of mutual distances between the Euclidean embedding of $\cT{W}$ is equal to $$\{\norm{\RM{w}_{s_1} - \RM{w}_{s_2}}_{\bT{R}^Q} = f^{\varepsilon}(\RM{w}_{s_1},\,\RM{w}_{s_2}): \RM{w}_{s_1},\,\,\RM{w}_{s_2} \in \cT{W}\},$$
where $f^{\varepsilon}$ is the symmetric proximity function obtained by replacing $\alpha$ with $\varepsilon$ in Eq.~\ref{eq:f_alpha}.
\end{proof}

%% file: intermediate.tex
In this section, we provide the intermediate results needed to find the value of $\varepsilon$ mentioned in Theorem~\ref{thm:main}.

\begin{definition}{(Connecting point)}
\label{def:cpoint}

Let us consider the set $\cT{X},\,\cT{Y}$ and a theoretical origin point $\RM{o}$ as characterized Notation~\ref{not:2sets}. Let us also consider the set of points $\cT{W} = \cT{X} \cup \cT{Y} \cup \{\RM{z}\}$ and $\cT{V}^{\RM{z}} = \cT{W} \cup \{\RM{z}\}$ as characterized in Notation~\ref{not:vz} and with associated proximity set $\cT{P}^{h^{\,\alpha}}_{\cT{V}^{\RM{z}}}$ (for $\alpha\geq0$) given in Eq.~\ref{eq:proxVz}, where $h^{\,\alpha}$ is given in Eq.~\ref{eq:h_alpha}.

As shown in Eq.~\ref{eq:h_alpha}, $h^{\,\alpha}$ depends on the choice of $\alpha$ through the proximity function $f^{\,\alpha}$ given by $$f^{\,\alpha}:\cT{V} \times \cT{V} \rightarrow \bT{R}^+,$$ satisfying Eq.~\ref{eq:f_alpha} and on the choice of proximity function $g$ that does not necessarily depends on $\alpha$, given by $$g:\cT{W} \times \{\RM{z}\} \rightarrow \bT{R}^+.$$

We say that $\RM{z} = \RM{z}^*$ is a ``connecting point'' associated with $h^{\,\alpha}$ if $\cT{V}^{\RM{z}^*}$, with respect to the proximity set $\cT{P}^{h^{\,\alpha}}_{\cT{V}^{\RM{z}^*}}$, can be embedded into an Euclidean space. 
\end{definition}

\begin{definition}{(Special proximity function for $\cT{V}^{\RM{z}}$)}
\label{def:tld_h}

Let us consider the sets $\cT{X},\,\cT{Y}$ with an associated proximity function $f$ and a theoretical origin point $\RM{o}$ (Notation~\ref{not:2sets}), and $\cT{W} = \cT{X}\cup \cT{Y} \cup \{\RM{o}\}$ (Notation~\ref{not:vz}).

We assume that one of the sets $\cT{X},\,\cT{Y}$ can be embedded into an Euclidean space and consider a point $\RM{w}_a$ belonging to that set. In all that follows, we restrict our work to the particular case $\RM{w}_a = \RM{x}_1 \in \cT{X}$ because it is always possible to rearrange the elements of the two sets $\cT{X}$ and $\cT{Y}$ and to relabel the sets to return to this particular case.

We denote
\begin{equation}
\label{eq:zeta_f}
    \zeta_f = \max_{\substack{l}} \quad \left\{\sum\limits_{\substack{s \\ s\neq l \\ \RM{w}_s \in \cT{X}, \RM{w}_l\in \cT{Y} \\ \text{ or } \\ \RM{w}_s \in \cT{Y}, \RM{w}_l\in \cT{X}}} \abs{\BT{m}_{l,s}^{\RM{x}_1}(\cT{W},\,f^{\,0})}\right\},
\end{equation}
where $f^{\,0}:\cT{W} \times \cT{W} \rightarrow \bT{R}^+,$ is the proximity function given in \ref{eq:f_alpha} when taking $\alpha = 0$ and $\BT{m}_{l,s}^{\RM{x}_1}(\cT{W},\,f^{\,0})$ is the $(l,s)$-component of the cosine law matrix $\BT{M}^{\RM{x}_1}(\cT{W},\,f^0)$ associated to $(\cT{W},\,f^{\,0})$ with reference point $\RM{x}_1$ and $f^{\,0}$ must be constructed such that $\zeta_f > 0$. In a nutshell, $\zeta_f$ is the maximum of the sum of the absolute values of the row elements of the anti-diagonal block matrices of $\BT{M}^{\RM{x}_1}(\cT{W},\,f^0)$.

For $\cT{V}^{\RM{z}}$ as characterized in Notation~\ref{not:vz}, we associate the special the proximity function $$h^{\,\alpha} = \tilde{h}^{\,\alpha}:\cT{V}^{\RM{z}} \times \cT{V}^{\RM{z}} \rightarrow \bT{R}^+,$$ given in Eq.~\ref{eq:h_alpha}, and in which $g$ is chosen as $$g = \tilde{g}_{\RM{x}_1} : \cT{W} \times \{\RM{z}\} \rightarrow \bT{R}^+,$$ such that for $s \in \{1, \ldots, \, S = M+N+1\}$ and $c_1, c_2 >0$
\begin{equation}
\label{eq:tld_g}
\tilde{g}_{\RM{x}_1}(\RM{w}_{s},\RM{z}) = 
\begin{cases}
(c_1\zeta_f)^{1/2} & \text{ if } \RM{w}_s = \RM{x}_1, \\

(c_2\zeta_f)^{1/2} & \text{ if } \RM{w_s} = \RM{o},\\

\left(d_{\cT{X}}(\RM{w}_{s},\RM{x}_1)^2 + c_2\zeta_f \right)^{1/2} & \text{ if } \RM{w}_s \neq \RM{x}_1 \text{ and } \RM{w}_s \in \cT{X},\\

\left(p^f_{\cT{X}\times \cT{Y}}(\RM{w}_{s},\RM{x}_1)^2 + c_2\zeta_f \right)^{1/2} & \text{ if } \RM{w}_s \in \cT{Y} ,\\
\end{cases}
\end{equation}
\end{definition}

\begin{notation}{(Split the cosine law matrix into a specific sum of matrices)}
\label{not:Msum}

Given $\alpha\geq 0$, $\zeta_f > 0$ and $c_1,\, c_2> 0$ and $\tilde{g}_{\RM{x}_1}$ as definite in Eq.~\ref{eq:tld_g} to construct the special proximity function $\tilde{h}^{\,\alpha}$ as in Def.~\ref{def:tld_h}. The cosine law matrix $\BT{M}^{\RM{z}}(\cT{V}^{\RM{z}},\,\tilde{h}^{\,\alpha})$ has the following components for $q_1,\, q_2 \in \{1, \ldots,\, Q = M+N+2\}$
$$\BT{m}_{q_1, q_2}^{\RM{z}}(\cT{V}^{\RM{z}},\,\tilde{h}^{\,\alpha}) = \dfrac{\tilde{h}^{\,\alpha}(\RM{v}_{q_1}, \RM{z})^2 + \tilde{h}^{\,\alpha}(\RM{v}_{q_2}, \RM{z})^2 - \tilde{h}^{\,\alpha}(\RM{v}_{q_1}, \RM{v}_{q_2})^2}{2}.$$
By definition, $\tilde{h}^{\,\alpha}$ depends on $\tilde{g}_{\RM{x}_1}$ (Eq.~\ref{eq:tld_g}), thus we obtain the following expression 

\begin{equation}
    \label{eq:details_mz}
    2\BT{m}_{q_1, q_2}^{\RM{z}}(\cT{V}^{\RM{z}},\,\tilde{h}^{\,\alpha}) = 
    \begin{cases}
        0                              
        & \hspace*{-5cm}\text{ if } \RM{v}_{q_1} = \RM{z} \text{ or } \RM{v}_{q_2} = \RM{z} \\

        c_1\zeta_f + c_2\zeta_f - \alpha 
        & \hspace*{-5cm}\text{ if } \RM{v}_{q_1}, \RM{v}_{q_2} \in \cT{X} \text{ and } (\RM{v}_{q_1} = \RM{x}_1, \RM{v}_{q_2}\neq \RM{x}_1) \text{ or } (\RM{v}_{q_1} \neq \RM{x}_1, \RM{v}_{q_2} = \RM{x}_1),\\

        2 c_1 \zeta_f 
        & \hspace*{-5cm}\text{ if } \RM{v}_{q_1} = \RM{x}_1, \RM{v}_{q_2} = \RM{x}_1,\\

        2c_2 \zeta_f - \alpha + \underbrace{d_{\cT{X}}(\RM{v}_{q_1}, \RM{x}_1)^2 + d_{\cT{X}}(\RM{v}_{q_2}, \RM{x}_1)^2 - d_{\cT{X}}(\RM{v}_{q_1}, \RM{v}_{q_2})^2}_{2\BT{m}^{\RM{x}_1}_{q_1, q_2}(\cT{X}, d_{\cT{X}})}
        & \text{ if } \RM{v}_{q_1}, \RM{v}_{q_2} \in \cT{X}\setminus\{\RM{x}_1\} \text{ and } \RM{v}_{q_1} \neq  \RM{v}_{q_2},\\

        2c_2\zeta_f - \alpha + d_{\cT{X}}(\RM{v}_{q_2}, \RM{x}_1)^2 - u_{\cT{X}}(\RM{v}_{q_2}, \RM{o}_{\cT{X}})^2
        & \text{ if } \RM{v}_{q_1} = \RM{o} \text{ and } \RM{v}_{q_2} \in \cT{X}\setminus \{x_1\}\\

        c_1\zeta_f + c_2\zeta_f - \alpha - u_{\cT{X}}(\RM{x}_1, \RM{o}_{\cT{X}})^2
        & \text{ if } \RM{v}_{q_1} = \RM{o} \text{ and } \RM{v}_{q_2} = \RM{x}_1,\\

        2c_2\zeta_f - \alpha + f(\RM{x}_1, \RM{v}_{q_2})^2 - u_{\cT{Y}}(\RM{v}_{q_2}, \RM{o}_{\cT{Y}})^2
        & \text{ if } \RM{v}_{q_1} = \RM{o} \text{ and } \RM{v}_{q_2} \in \cT{Y},\\

        2c_2\zeta_f - \alpha + d_{\cT{X}}(\RM{v}_{q_1}, \RM{x}_1)^2 - u_{\cT{X}}(\RM{v}_{q_1}, \RM{o}_{\cT{X}})^2
        & \text{ if } \RM{v}_{q_1} \in \cT{X}\setminus \{\RM{x}_1\} \text{ and } \RM{v}_{q_2} = \RM{o},\\

        c_1\zeta_f + c_2\zeta_f - \alpha - u_{\cT{X}}(\RM{x}_{1}, \RM{o}_{\cT{X}})^2
        & \text{ if } \RM{v}_{q_1} = \RM{x}_1\text{ and } \RM{v}_{q_2} = \RM{o},\\

        2c_2\zeta_f - \alpha + f(\RM{x}_1,\RM{v}_{q_1})^2 - u_{\cT{Y}}(\RM{v}_{q_1}, \RM{o}_{\cT{Y}})^2
        & \text{ if } \RM{v}_{q_1} \in \cT{Y} \text{ and } \RM{v}_{q_2} = \RM{o},\\

        2c_2\zeta_f - \alpha + f(\RM{x}_1, \RM{v}_{q_1})^2 + f(\RM{x}_1, \RM{v}_{q_2})^2 - d_{\cT{Y}}(\RM{v}_{q_1}, \RM{v}_{q_2})^2
        & \text{ if } \RM{v}_{q_1} \text{ and } \RM{v}_{q_2} \in \cT{Y}\text{ and } \RM{v}_{q_1} \neq  \RM{v}_{q_2},\\

        c_1\zeta_f + c_2\zeta_f - \alpha 
        & \text{ if } \RM{v}_{q_1} = \RM{x}_1 \text{ and } \RM{v}_{q_2} \in \cT{Y},\\

        2c_2\zeta_f - \alpha  + \underbrace{d_{\cT{X}}(\RM{v}_{q_1}, \RM{x}_1)^2 + f(\RM{x}_1,\RM{v}_{q_2})^2 - f(\RM{v}_{q_1}, \RM{v}_{q_2})^2}_{2\BT{m}^{\RM{x}_1}_{q_1,q_2}(\cT{W}, f^0)}
        & \text{ if } \RM{v}_{q_1} \in \cT{X}\setminus \{\RM{x}_1\}  \text{ and } \RM{v}_{q_2} \in \cT{Y},\\

        c_1\zeta_f + c_2\zeta_f - \alpha 
        & \text{ if } \RM{v}_{q_1} \in \cT{Y} \text{ and } \RM{v}_{q_2} = \RM{x}_1,\\

        2c_2\zeta_f - \alpha + \underbrace{d_{\cT{X}}(\RM{v}_{q_2}, \RM{x}_1)^2 + f(\RM{x}_1, \RM{v}_{q_1})^2 - f(\RM{v}_{q_2}, \RM{v}_{q_1})^2}_{2\BT{m}^{\RM{x}_1}_{q_1,q_2}(\cT{W}, f^0)}
        & \text{ if } \RM{v}_{q_1} \in \cT{Y}  \text{ and } \RM{v}_{q_2} \in \cT{X}\setminus \{\RM{x}_1\} ,\\

        2c_2\zeta_f + \underbrace{2d_{\cT{X}}(\RM{v}_{q_1}, \RM{x}_1)^2}_{2\BT{m}^{\RM{x}_1}_{q_1, q_1}(\cT{X}, d_{\cT{X}})}
        & \text{ if } \RM{v}_{q_1} = \RM{v}_{q_2} \text{ and } \RM{v}_{q_1},\RM{v}_{q_2}\in \cT{X}\setminus\{\RM{x}_1\},\\
        2c_2\zeta_f + 2f(\RM{x}_1,\RM{v}_{q_1})^2
        & \text{ if } \RM{v}_{q_1} = \RM{v}_{q_2} \text{ and } \RM{v}_{q_1},\RM{v}_{q_2}\in \cT{Y},
    \end{cases}
\end{equation}

Then, to support our proofs, we split $\BT{M}^{\RM{z}}(\cT{V}^{\RM{z}},\,\tilde{h}^{\,\alpha})$ in a specific way as the follows
\begin{equation}
\label{eq:Msum}
    \BT{M}^{\RM{z}}(\cT{V}^{\RM{z}},\,\tilde{h}^{\,\alpha}) = \BT{M}^{\RM{x}_1}_{\cT{X}} \, + \, \BT{G}^{\RM{x}_1}_{\cT{X}\times \cT{Y}}(\alpha) \, + \,  \BT{C}^{\RM{x}_1}_{\cT{X} \times \{\RM{o}_{\cT{X}}\}}(\alpha) \, + \, \BT{C}^{\RM{x}_1}_{\cT{Y} \times \{\RM{o}_{\cT{Y}}\}}(\alpha),
\end{equation}
where
\begin{equation}
\label{eq:Mx1}
    \BT{M}^{\RM{x}_1}_{\cT{X}} =
    \begin{pmatrix}
         0 
         & \BT{0}_{1\times M} 
         & 0
         & \BT{0}_{1 \times N}\\
         
         \BT{0}_{M\times 1}
         & \red{\BT{M}^{\RM{x}_1}(\cT{X}, d_{\cT{X}})}
         & \BT{0}_{M\times 1}
         & \BT{0}_{M\times N}\\
         
         0 
         & \BT{0}_{1\times M}
         & 0
         & \BT{0}_{1\times N}\\
         
         \BT{0}_{N\times 1}
         & \BT{0}_{N\times M}
         & \BT{0}_{N\times 1}
         & \BT{0}_{N\times N}
    \end{pmatrix}
\end{equation}

where $\bullet \times \bullet$ indicate the dimensions on the matrices, $\BT{0}_{\bullet \times \bullet}$ are null martices, and $\BT{M}^{\RM{x}_1}(\cT{X}, d_{\cT{X}})$ is the cosine law matrix associated with $(\cT{X}, d_{\cT{X}})$ and with reference point $\RM{x}_1$. Next,

\begin{equation}
\label{eq:GXY}
    \BT{G}^{\RM{x}_1}_{\cT{X}\times\cT{Y}}(\alpha) =
    \begin{pmatrix}
             0 
             & \BT{0}_{1\times M} 
             & 0
             & \BT{0}_{1 \times N}\\
             
             \BT{0}_{M\times 1}
             & \red{\dfrac{c_2\zeta_f}{2} \BT{I}_{M\times M}}
             & \BT{0}_{M\times 1}
             & \red{\BT{B}^{\RM{x}_1}_{\cT{X}\times \cT{Y}}(\alpha)}\\
             0 
             & \BT{0}_{1\times M}
             & 0
             & \BT{0}_{1\times N}\\
             
             \BT{0}_{N\times 1}
             & \red{\left(\BT{B}^{\RM{x}_1}_{\cT{X}\times \cT{Y}}(\alpha)\right)^{\intercal}}
             & \BT{0}_{N\times 1}
             & \red{\dfrac{c_2\zeta_f}{2} \BT{I}_{N\times N}},
        \end{pmatrix}  
\end{equation}
where $I_{\bullet \times \bullet}$ is the identity matrix and the submatrix $\BT{B}^{\RM{x}_1}_{\cT{X}\times\cT{Y}}(\alpha)$ is a $M\times N$ matrix with components
\begin{equation}
\label{eq:fXYa}
 2\,\BT{b}^{\RM{x}_1}_{\cT{X} \times \cT{Y}}(\alpha)_{m, n} =
 \begin{cases}
  c_1\zeta_f + c_2\zeta_f - \alpha  
  & \text{ if } \RM{x}_m = \RM{x}_1,\\
  
  2c_2\zeta_f - \alpha + \underbrace{d_{\cT{X}}(\RM{x}_m, \RM{x}_1)^2+ f(\RM{x}_1, \RM{y}_n)^2  - f(\RM{x}_m,\RM{y}_n)^2}_{2\BT{m}^{\RM{x}_1}_{m,n+M+2}(\cT{W}, f^0)} 
  & \text{ if } \RM{x}_m \neq \RM{x}_1.\\
 \end{cases}  
\end{equation}
Next, 
\begin{equation}
\label{eq:CXo}
    \BT{C}^{\RM{x}_1}_{\cT{X} \times \{\RM{o}_{\cT{X}}\}}(\alpha) =
      \begin{pmatrix}
         0 
         & \BT{0}_{1\times M} 
         & 0
         & \BT{0}_{1 \times N}\\
         
         \BT{0}_{M\times 1}
         & \red{\BT{A}_{\cT{X}}(\alpha) - \dfrac{c_2\zeta_f}{2} \BT{I}_{M\times M}}
         & \red{\left(\BT{C}^{\RM{x}_1}_{\RM{o}_{\cT{X}}}(\alpha)\right)^{\,\intercal}} 
         & \BT{0}_{M\times N}\\
         
         0 
         & \red{\BT{C}^{\RM{x}_1}_{\RM{o}_{\cT{X}}}(\alpha)}
         & \red{c_2\zeta_f/2}
         & \BT{0}_{1\times N}\\
         
         \BT{0}_{N\times 1}
         & \BT{0}_{N\times M}
         & \BT{0}_{N\times 1}
         & \BT{0}_{N\times N}
    \end{pmatrix},  
\end{equation}
where $\BT{C}^{\RM{x}_1}_{\RM{o}_{\cT{X}}}(\alpha)$ is an $1\times M$ matrices with components
\begin{equation}
\label{eq:v_1oX}
    2\,\BT{c}^{\RM{x}_1}_{\RM{o}_{\cT{X}}}(\alpha)_{1,m} = 
    \begin{cases}
    c_1\zeta_f + c_2\zeta_f - \alpha - u_\cT{X}(\RM{x}_1,\, \RM{o}_{\cT{X}})^2, 
    & \text{ if } \RM{x}_m = \RM{x}_1,\\
    
      2\,c_2\zeta_f  - \alpha + d_{\cT{X}}(\RM{x}_m, \RM{x}_1)^2 - u_\cT{X}(\RM{x}_{m},\, \RM{o}_{\cT{X}})^2, 
    & \text{ if } \RM{x}_m \neq \RM{x}_1,\\
    \end{cases}
\end{equation}
and 
\begin{equation}
\label{eq:A_X}
    2\,\BT{a}_{\cT{X}}(\alpha)_{m_1, m_2} =
    \begin{cases}
     2 c_1\zeta_f & \text{ if } \RM{x}_{m_1} = \RM{x}_{m_2} = \RM{x}_1\\
     
       c_1\zeta_f + c_2\zeta_f - \alpha
       & \text{ if } \RM{x}_{m_1} = \RM{x}_1, \, \RM{x}_{m_1} \neq \RM{x}_{m_2} \text{ or vice-versa},\\ 
       
        2c_2\zeta_f - \alpha
        &\text{ if } \RM{x}_{m_1},\, \RM{x}_{m_2} \in \cT{X}\setminus\{\RM{x}_1\} \text{ and } \RM{x}_{m_1}\neq \RM{x}_{m_2},\\
      2c_2\zeta_f
        & \text{ if } \RM{x}_{m_1} = \RM{x}_{m_2}  \text{ and } \RM{x}_{m_1},\, \RM{x}_{m_2} \in \cT{X}\setminus\{\RM{x}_1\},
    \end{cases}
\end{equation}
Lastly,
\begin{equation}
\label{eq:CYo}
    \BT{C}^{\RM{x}_1}_{\cT{Y} \times \{\RM{o}_{\cT{Y}}\}}(\alpha) =
      \begin{pmatrix}
         0 
         & \BT{0}_{1\times M} 
         & 0
         & \BT{0}_{1 \times N}\\
         
         \BT{0}_{M\times 1}
         & \BT{0}_{M\times M}
         & \BT{0}_{M\times 1}
         & \BT{0}_{M\times N}\\
         
         0 
         & \BT{0}_{1\times M}
         & \red{c_2\zeta_f/2}
         & \red{\BT{C}^{\RM{x}_1}_{\RM{o}_{\cT{Y}}}(\alpha)}\\
         
         \BT{0}_{N\times 1}
         & \BT{0}_{N\times M}
         & \red{\left(\BT{C}^{\RM{x}_1}_{\RM{o}_{\cT{Y}}}(\alpha)\right)^{\,\intercal}}
         & \red{\BT{A}^{\RM{x}_1}_{\cT{Y}}(\alpha) - \dfrac{c_2\zeta_f}{2} \BT{I}_{N\times N}}
    \end{pmatrix},  
\end{equation}
where $\BT{C}^{\RM{x}_1}_{\RM{o}_{\cT{Y}}}(\alpha)$ is an $1 \times N$ matrix with components
\begin{equation}
\label{eq:v_1oY}
    2\,\BT{c}^{\RM{x}_1}_{\RM{o}_{\cT{Y}}}(\alpha)_{1,n} = 
    2\,c_2\zeta_f - \alpha + f(\RM{x}_1,\, \RM{y}_n)^2 - u_\cT{Y}(\RM{y}_{n},\, \RM{o}_{\cT{Y}})^2.
\end{equation}
\end{notation}

and $\BT{A}^{\RM{x}_1}_{\cT{Y}}(\alpha)$ is a $N \times N$ matrix with components

\begin{equation}
\label{eq:A_x1Y}
 2\,\BT{a}^{\RM{x}_1}_{\cT{Y}}(\alpha)_{n_1, n_2}
    =
    \begin{cases}
    2c_2\zeta_f + 2f(\RM{x}_1, \RM{y}_{n_1})^2
    & \text{ if } \RM{y}_{n_1} = \RM{y}_{n_2},\\

    2c_2\zeta_f - \alpha + f(\RM{x}_1,\, \RM{y}_{n_1})^2 + f(\RM{x}_1,\, \RM{y}_{n_2})^2 - d_{\cT{Y}}(\RM{y}_{n_1},\, \RM{y}_{n_2})^2, 
    & \text{ otherwise }
    \end{cases}
\end{equation}

\begin{theorem}
\label{thm:vareps}

Let us consider the sets $\cT{X},\,\cT{Y}$ with an associated proximity function $f$ and a theoretical origin point $\RM{o}$ (Notation~\ref{not:2sets}), $\cT{W} = \cT{X}\cup \cT{Y} \cup \{\RM{o}\}$ with $S$ elements and $\cT{V}^{\RM{z}} = \cT{W}\cup \{\RM{z}\}$ with $Q$ elements as characterized in Notation~\ref{not:vz}.

For $\alpha> 0$, $\zeta_f > 0$ given in Eq.~\ref{eq:zeta_f} and $c_1,\, c_2> 0$ defining $\tilde{g}_{\RM{x}_1}$ in Eq.~\ref{eq:tld_g}, we construct the special proximity function $\tilde{h}^{\,\alpha}$ associated with a theoretical point $\RM{z}$ in Def.~\ref{def:tld_h}.

If we choose a third number $c_3 > 0$ such that
\begin{itemize}
    \item[(i)] for some non-negative constants $K$ and $K'$ dependent on $c_1$ and $c_2$, $c_1,\, c_2,\, c_3>0$ are solutions of 
    \begin{equation}
        \label{eq:cond_i_c123}
        \begin{cases}
        c_1 - \max(M,N) \abs{2 + c_1+ c_2 - c_3} - K(c_1, c_2) \geq 0,\\
        c_2  - \max(M, N) \abs{2 + c_1 + c_2 - c_3} - K'(c_1, c_2) \geq 0.\\
        \end{cases}
\end{equation}  
\end{itemize}
Then, if we choose $\alpha = \varepsilon = c_3\zeta_f$, the theoretical point $\RM{z} = \RM{z^*}$ associated with $\tilde{h}^{\,\varepsilon}$ is a connecting point.
\end{theorem}

\begin{remark}{(Solution of condition (i) Theorem~\ref{thm:vareps})}

Besides the equations listed in conditions (i) of Theorem~\ref{thm:vareps}, there might still be other choices of $c_1, c_2, c_3>0$ satisfying the condition for $\varepsilon = c_3\zeta_f$ because the conditions are sufficient but not necessary for $\BT{M}^{\RM{z}^*}(\cT{V}^{\RM{z}^*},\,\tilde{h}^{\,\varepsilon})$ to be positive semi-definite. For example, the triplets
\begin{equation}
\label{eq:eg}
    \left(
        c_1 = 1,\,
        c_2 = 1/2,\,
        c_3 = \dfrac{2c_1 + c_2 - b}{a}\right) \text{ with } a = 1 - \dfrac{1}{M+N} \text{ and } b = \dfrac{2c_2}{M+N}.
\end{equation}
was used to successfully find qualitative Euclidean embedding for specific set examples (see Code availability). However, the triplet in Eq.~\ref{eq:eg} does not satisfy the conditions in (i) for $(M,N) = (5, 4)$ or $(M, N) = (3, 5)$. This triplet was found in previous trials and errors. 

In general, it is evident that the following triplet is a general solution to condition (i)
\begin{equation}
\label{eq:eg2}
    \left(
        c_1 \geq K(c_1, c_2),\,
        c_2 \geq K'(c_1, c_2),\,
        c_3 = 2 + c_1 + c_2
        \right).
\end{equation}
which is a section of the upper half-plane. In the proof of Theorem~\ref{thm:vareps}, we show that $K$ and $K'$ also depends on distance/proximity functions $d_{\cT{X}},\, d_{\cT{Y}},$ and $f$ which are all finite-valued functions. Thus, for any particular problem, there exists a large enough value $\Lambda > 0 $ such that the triplet 
\begin{equation}
\label{eq:eg3}
    \left(
        c_1 = \Lambda,\,
        c_2 = 2 \Lambda,\,
        c_3 = 2 + c_1 + c_2.
        \right),
\end{equation}
is a general solution of condition (i). Therefore the following algorithm converges to a particular solution.
\begin{algorithm}
\caption{An algorithm that finds a solution for condition (i) Theorem~\ref{thm:vareps}} \label{alg:c1c2c3}
$c_1 \gets $ random $>$ 0\;
$c_2 \gets $ random $>$ 0\;
$c_3 \gets $ random $>$ 0\;
\While{condition (i) FALSE}{
    $c_1 \gets c_2$\;
    $c_2 \gets 2c_1$\;
    $c_3 \gets 2 + c_1 + c_2$\;
}
\end{algorithm}

In the future, the above algorithm will be used for general problems.
\end{remark}

\begin{proof}[proof of Theorem~\ref{thm:vareps}]
It is sufficient to find the conditions which the cosine law matrix $\BT{M}^{\RM{z}}(\cT{V}^{\RM{z}},\,\tilde{h}^{\,\alpha})$ is for positive semi-definiteness. If the matrices on the right-hand side of Eq.~\ref{eq:Msum} are positive semi-definite, then the cosine law matrix $\BT{M}^{\RM{z}}(\cT{V}^{\RM{z}},\,\tilde{h}^{\,\alpha})$ is positive semi-definite (Lemma~\ref{lm:sumPSD}). It remains the conditions in (i) of Theorem~\ref{thm:vareps} are satisfied, then for a specific construction of $\RM{z} = \RM{z}^*$ obtained by a specific choice of $\alpha = \varepsilon = c_3\zeta_f$,  $\BT{M}^{\RM{x}_1}_{\cT{X}}$ (Eq.~\ref{eq:Mx1}),  $\BT{G}^{\RM{x}_1}_{\cT{X}\times \cT{Y}}(\varepsilon)$ (Eq.~\ref{eq:GXY}),  $\BT{C}^{\RM{x}_1}_{\cT{X} \times \{\RM{o}_{\cT{X}}\}}(\varepsilon)$ (Eq.~\ref{eq:v_1oX})  $\BT{C}^{\RM{x}_1}_{\cT{Y} \times \{\RM{o}_{\cT{Y}}\}}(\varepsilon)$ (Eq.~\ref{eq:v_1oY}) are positive semi-definite.

Firstly, following condition (i) of Theorem~\ref{thm:main}, the distance set $\cT{D}_{\cT{X}}$, associated with $d_{\cT{X}}$, is a set of mutual Euclidean distances between the elements of $\cT{X}$, the cosine law matrix $\BT{M}^{\RM{x}_1}(\cT{X}, d_{\cT{X}})$ is positive semi-definite (if is $\cT{D}_{\cT{Y}}$ that is a mutual Euclidean distance, then we swap the labels of the two sets to recover the same statement). Therefore  $\BT{M}^{\RM{x}_1}_{\cT{X}}$ is also positive semi-definite because it only has $\BT{M}^{\RM{x}_1}(\cT{X}, d_{\cT{X}})$ as non-zero submatrix. 

Secondly, for $\BT{G}^{\RM{x}_1}_{\cT{X}\times\cT{Y}}(\varepsilon)$, the rows having all 0 components can also be ignored since they correspond to the eigenvalue 0. Let $\lambda$ be an eigenvalue of $\BT{G}^{\RM{x}_1}_{\cT{X}\times\cT{Y}}(\varepsilon)$ corresponding to a non-zero row. The Ger\v{s}gorin circle theorem (Theorem~\ref{thm:G}) states that
\begin{equation}
\label{eq:iLbdG}
\exists\: q_{\lambda} \in \{1, \ldots, Q\}, \text{ such that } \abs{\lambda - \BT{g}^{\RM{x}_1}_{\cT{X}\times\cT{Y}}(\varepsilon)_{q_{\lambda}, q_{\lambda}}} \,\leq\, R_{q_{\lambda}}(\BT{G}^{\RM{x}_1}_{\cT{X}\times\cT{Y}}(\varepsilon)) 
\end{equation}
where 
\begin{equation*}
    R_{q_{\lambda}}(\BT{G}^{\RM{x}_1}_{\cT{X}\times\cT{Y}}(\varepsilon)) = 
    \sum\limits_{j\neq q_{\lambda}}\abs{\BT{g}^{\RM{x}_1}_{\cT{X}\times\cT{Y}}(\varepsilon)_{q_{\lambda},j}},
\end{equation*}
which implies
\begin{equation*}
 R_{q_{\lambda}}(\BT{G}^{\RM{x}_1}_{\cT{X}\times\cT{Y}}(\varepsilon))    =
        \begin{cases}
         \sum\limits_{j = M+3}^{Q} 
         \abs{\BT{b}^{\RM{x}_1}_{\cT{X}\times \cT{Y}}(\varepsilon)_{q_{\lambda}, j-{M+2}}}
         & \text{ if } q_{\lambda} = {2, \ldots,\, M+1},\\
         \\
         \sum\limits_{j = 2}^{M+1}
         \abs{\BT{b}^{\RM{x}_1}_{\cT{X}\times \cT{Y}}(\varepsilon)_{j, q_{\lambda} - (M+2)}}
         & \text{ if } q_{\lambda} = {M+3, \ldots,\, Q=M+N+2}.
        \end{cases}
\end{equation*}
According to Eq.~\ref{eq:zeta_f}~and~\ref{eq:fXYa}, $\zeta_f>0$, and $\zeta_f$ and $\BT{b}^{\RM{x}_1}_{\cT{X},\cT{Y}}(\varepsilon)_{\bullet \times \bullet}$ are related. Then, by the triangle inequality, we obtain the following expression

\begin{equation*}
        \begin{cases}
        R_{q_{\lambda}}(\BT{G}^{\RM{x}_1}_{\cT{X}\times\cT{Y}}(\varepsilon)) \leq
         \dfrac{N}{2}\left(\abs{c_1\zeta_f + c_2\zeta_f - \varepsilon}\right)
         & \text{ if } \RM{v}_{q_{\lambda}} = \RM{x}_1 \text{ and } q_{\lambda} = 2, \ldots,\, M+1,\\
         \\
        
        R_{q_{\lambda}}(\BT{G}^{\RM{x}_1}_{\cT{X}\times\cT{Y}}(\varepsilon)) \leq
         \dfrac{N}{2}\abs{2c_2\zeta_f - \varepsilon} + \dfrac{1}{2}\zeta_f
         & \text{ if } \RM{v}_{q_{\lambda}} \neq \RM{x}_1 \text{ and } q_{\lambda} = 2, \ldots,\, M+1,\\
         \\
         R_{q_{\lambda}}(\BT{G}^{\RM{x}_1}_{\cT{X}\times\cT{Y}}(\varepsilon)) \leq 
         
         \dfrac{M}{2}\abs{2c_2\zeta_f - \varepsilon} + \dfrac{1}{2}\zeta_f
         & \text{ if } q_{\lambda} = {M+3, \ldots,\, Q=M+N+2}.
        \end{cases}
\end{equation*}
Combining the above equation with Eq.~\ref{eq:iLbdG} gives
\begin{equation*}
\label{eq:GBound_0}
        \begin{cases}
          \abs{\lambda - \BT{g}^{\RM{x}_1}_{\cT{X}\times\cT{Y}}(\varepsilon)_{q_{\lambda}, q_{\lambda}}} \leq
         \dfrac{N}{2} \abs{c_1\zeta_f + c_2\zeta_f - \varepsilon}
         & \text{ if } \RM{v}_{q_{\lambda}} = \RM{x}_1 \text{ and } q_{\lambda} = 2, \ldots,\, M+1,\\
         \\
        
         \abs{\lambda - \BT{g}^{\RM{x}_1}_{\cT{X}\times\cT{Y}}(\varepsilon)_{q_{\lambda}, q_{\lambda}}} \leq
         \dfrac{N}{2} \abs{2c_2\zeta_f - \varepsilon} 
         + 
        \dfrac{1}{2}\zeta_f
         & \text{ if } \RM{v}_{q_{\lambda}} \neq \RM{x}_1 \text{ and } q_{\lambda} = 2, \ldots,\, M+1,\\
         \\
          \abs{\lambda - \BT{g}^{\RM{x}_1}_{\cT{X}\times\cT{Y}}(\varepsilon)_{q_{\lambda}, q_{\lambda}}} \leq 
         
         \dfrac{M}{2}\abs{2c_2\zeta_f - \varepsilon}
         + 
         \dfrac{1}{2}\zeta_f
         & \text{ if } q_{\lambda} = {M+3, \ldots,\, Q=M+N+2}.
        \end{cases}
\end{equation*}
What is more important for us here is to derive a lower bound for $\lambda$ from the above equation. We, thus, have the following
\begin{equation*}
\label{eq:GBound_1}
        \begin{cases}
          \lambda \geq \BT{g}^{\RM{x}_1}_{\cT{X}\times\cT{Y}}(\varepsilon)_{q_{\lambda}, q_{\lambda}} 
          - \dfrac{N}{2} \abs{c_1\zeta_f + c_2\zeta_f - \varepsilon}
         & \text{ if } \RM{v}_{q_{\lambda}} = \RM{x}_1 \text{ and } q_{\lambda} = 2, \ldots,\, M+1,\\
         \\
        
        \lambda \geq \BT{g}^{\RM{x}_1}_{\cT{X}\times\cT{Y}}(\varepsilon)_{q_{\lambda}, q_{\lambda}}
         - \dfrac{N}{2} \abs{2c_2\zeta_f - \varepsilon} 
         - \dfrac{1}{2}\zeta_f
         & \text{ if } \RM{v}_{q_{\lambda}} \neq \RM{x}_1 \text{ and } q_{\lambda} = 2, \ldots,\, M+1,\\
         \\
          \lambda \geq \BT{g}^{\RM{x}_1}_{\cT{X}\times\cT{Y}}(\varepsilon)_{q_{\lambda}, q_{\lambda}} 
         
         - \dfrac{M}{2}\abs{2c_2\zeta_f - \varepsilon}
         - \dfrac{1}{2}\zeta_f
         & \text{ if } q_{\lambda} = {M+3, \ldots,\, Q=M+N+2},
        \end{cases}
\end{equation*}
Now, we set $\varepsilon = c_3\zeta_f$ and replace $\BT{g}^{\RM{x}_1}_{\cT{X}\times\cT{Y}}(\varepsilon)_{q_{\lambda}, q_{\lambda}}$ by its value to obtain 
\begin{equation*}
\label{eq:GBound_2}
        \begin{cases}
          \lambda \geq \dfrac{c_1\zeta_f}{2} 
          - \dfrac{N}{2} \abs{c_1\zeta_f + c_2\zeta_f - c_3\zeta_f}
         & \text{ if } \RM{v}_{q_{\lambda}} = \RM{x}_1 \text{ and } q_{\lambda} = 2, \ldots,\, M+1,\\
         \\
        
         \lambda \geq \dfrac{c_2\zeta_f}{2} 
         - \dfrac{N}{2} \abs{2c_2\zeta_f - c_3\zeta_f} 
         - \dfrac{1}{2}\zeta_f
         & \text{ if } \RM{v}_{q_{\lambda}} \neq \RM{x}_1 \text{ and } q_{\lambda} = 2, \ldots,\, M+1,\\
         \\
         \lambda \geq \dfrac{c_2\zeta_f}{2}
         - \dfrac{M}{2}\abs{2c_2\zeta_f - c_3\zeta_f}
         - \dfrac{1}{2}\zeta_f
         & \text{ if } q_{\lambda} = {M+3, \ldots,\, Q=M+N+2}.
        \end{cases}
\end{equation*}
Using the triangle inequality, and the change of expression $c_1\zeta_f + c_2\zeta_f - \varepsilon =  2\zeta_f + c_1\zeta_f + c_2\zeta_f - \varepsilon - 2\zeta_f $ and $2c_2\zeta_f - \varepsilon = 2\zeta_f + c_1\zeta_f + c_2\zeta_f - \varepsilon + c_2\zeta_f - c_1\zeta_f - 2\zeta_f$ we have
\begin{equation*}
\label{eq:GBound_2a}
        \begin{cases}
         \dfrac{c_1\zeta_f}{2} 
          - \dfrac{N }{2} \left(\abs{2\zeta_f + c_1\zeta_f + c_2\zeta_f - \varepsilon} + 2\zeta_f\right) \geq 0
         \\
         \\
         \dfrac{c_2\zeta_f}{2} 
         - \dfrac{N}{2} \left(\abs{2\zeta_f + c_1 + c_2 - c_3} + \abs{c_2\zeta_f - c_1\zeta_f + 2\zeta_f}\right)
         - \dfrac{1}{2} \geq 0
         \\
         \\
          \dfrac{c_2\zeta_f}{2} 
         - \dfrac{M}{2}\left(\abs{2\zeta_f + c_1\zeta_f + c_2\zeta_f - \varepsilon} + \abs{c_2\zeta_f - c_1\zeta_f + 2\zeta_f}\right)
         - \dfrac{1}{2} \geq 0
        \end{cases}
\end{equation*}
and thus the conditions for which $\BT{G}^{\RM{x}_1}_{\cT{X}\times \cT{Y}}(\varepsilon=c_3\zeta_f)$ is positive semi-definite are
\begin{equation}
\label{eq:GBound}
        \begin{cases}
         c_1 - N \abs{2 + c_1 + c_2 - c_3} - K_0(c_1, c_2) \geq 0\\
         c_2 - \max(M, N) \abs{2 + c_1 + c_2 - c_3} - K_1(c_1, c_2) \geq 0
        \end{cases}
\end{equation}
for some non-negative constants $K_0(c_1, c_2)$ and $K_1(c_1, c_2)$ dependent on $c_1$ and $c_2$. And this condition is summarized in condition (i) of Theorem~\ref{thm:vareps}. 

Lastly, we investigate the conditions for which $\BT{C}^{\RM{x}_1}_{\cT{X} \times \{\RM{o}_{\cT{X}}\}}(\varepsilon)$ and $\BT{C}^{\RM{x}_1}_{\cT{Y} \times \{\RM{o}_{\cT{Y}}\}}(\varepsilon)$ are positive semi-definite. Let $\gamma$ and $\delta$ be eigenvalues corresponding to the non-zero rows $\BT{C}^{\RM{x}_1}_{\cT{X} \times \{\RM{o}_{\cT{X}}\}}(\varepsilon)$ and $\BT{C}^{\RM{x}_1}_{\cT{Y} \times \{\RM{o}_{\cT{Y}}\}}(\varepsilon)$, respectively. Then, by the Ger\v{s}gorin circle theorem (Theorem~\ref{thm:G}), we have the following

\begin{equation}
\label{eq:gamma_delta}
\begin{cases}
    \exists\, q_{\gamma}, \in \{1, \ldots, Q\}, \text{ such that } \abs{\gamma - \BT{c}^{\RM{x}_1}_{\cT{X} \times \{\RM{o}_{\cT{X}}\}}(\varepsilon)_{q_{\gamma}, q_{\gamma}}} \,\leq\, R_{q_{\gamma}}(\BT{C}^{\RM{x}_1}_{\cT{X} \times \{\RM{o}_{\cT{X}}\}}(\varepsilon)),\\
    
    \exists\, q_{\delta}, \in \{1, \ldots, Q\}, \text{ such that } \abs{\delta - \BT{c}^{\RM{x}_1}_{\cT{Y} \times \{\RM{o}_{\cT{Y}}\}}(\varepsilon)_{q_{\delta}, q_{\delta}}} \,\leq\, R_{q_{\delta}}(\BT{C}^{\RM{x}_1}_{\cT{Y} \times \{\RM{o}_{\cT{Y}}\}}(\varepsilon)),
\end{cases}
\end{equation}
where
\begin{equation}
\begin{cases}
    R_{q_{\gamma}}(\BT{C}^{\RM{x}_1}_{\cT{X} \times \{\RM{o}_{\cT{X}}\}}(\varepsilon)) 
    = 
    \sum\limits_{j\neq q_{\gamma}} \abs{\BT{c}^{\RM{x}_1}_{\cT{X} \times \{\RM{o}_{\cT{X}}\}}(\varepsilon)_{q_{\gamma}, j}},\\
    R_{q_{\delta}}(\BT{C}^{\RM{x}_1}_{\cT{Y} \times \{\RM{o}_{\cT{Y}}\}}(\varepsilon))
    = \sum\limits_{j\neq q_{\delta}} \abs{\BT{c}^{\RM{x}_1}_{\cT{Y} \times \{\RM{o}_{\cT{Y}}\}}(\varepsilon)_{q_{\delta}, j}}.
\end{cases}
\end{equation}
By replacing the components with their respective values, we get
\begin{equation*}
    \begin{cases}
    R_{q_{\gamma}}(\BT{C}^{\RM{x}_1}_{\cT{X} \times \{\RM{o}_{\cT{X}}\}}(\varepsilon)) 
    = 
    \sum\limits_{\substack{j = 2\\j\neq q_{\gamma}}}^{M+1} \abs{\BT{a}_{\cT{X}}(\varepsilon)_{q_{\gamma}, j}} 
    +
    \frac{1}{2} \abs{c_1\zeta_f + c_2\zeta_f - \varepsilon - u_{\cT{X}}(\RM{x}_1, \RM{o}_{\cT{X}})^2},
    & \hspace*{-2.5cm} \text{ if } \RM{v}_{q_{\gamma}} = \RM{x}_1 \text{ i.e., } q_{\gamma} = 2,
    \\
    \\
    R_{q_{\gamma}}(\BT{C}^{\RM{x}_1}_{\cT{X} \times \{\RM{o}_{\cT{X}}\}}(\varepsilon)) 
    = 
    \sum\limits_{\substack{j = 2\\j\neq q_{\gamma}}}^{M+1} \abs{\BT{a}_{\cT{X}}(\varepsilon)_{q_{\gamma}, j}} 
    +
    \frac{1}{2}  \abs{2c_2\zeta_f - \varepsilon + d_{\cT{X}}(\RM{v}_{q_{\gamma}}, \RM{x}_1)^2 - u_{\cT{X}}(\RM{x}_1, \RM{o}_{\cT{X}})^2},
    & \hspace*{-2.5cm}\text{ if } q_{\gamma} = 3, \ldots, \, M+1,
    \\
    \\
     R_{q_{\gamma}}(\BT{C}^{\RM{x}_1}_{\cT{X} \times \{\RM{o}_{\cT{X}}\}}(\varepsilon)) 
    = \frac{1}{2}  \sum\limits_{m=1}^{M}\abs{2c_2\zeta_f - \varepsilon + d_{\cT{X}}(\RM{v}_{q_{\gamma}}, \RM{x}_1)^2 - u_{\cT{X}}(\RM{x}_1, \RM{o}_{\cT{X}})^2},
    & \hspace*{-3.5cm}\text{ if } q_{\gamma} = M+2,
    \\
    \\
    R_{q_{\delta}}(\BT{C}^{\RM{x}_1}_{\cT{Y} \times \{\RM{o}_{\cT{Y}}\}}(\varepsilon))
    = 
    \sum\limits_{\substack{j = M+3\\j\neq q_{\delta}}}^{Q}
         \abs{\BT{a}^{\RM{x}_1}_{\cT{Y}}(\varepsilon)_{q_{\delta} - (M+2), j - (M+2)}}
    +
    \frac{1}{2} \abs{2c_2\zeta_f - \varepsilon + f(\RM{x}_1, \RM{v}_{q_{\delta}})^2- u_{\cT{Y}}(\RM{v}_{q_{\delta}}, \RM{o}_{\cT{Y}})^2},
    \\
    & \hspace*{-5cm}\text{ if } q_{\delta} = M+3, \ldots,\,Q = M+N+2,
    \\
    \\
    R_{q_{\delta}}(\BT{C}^{\RM{x}_1}_{\cT{Y} \times \{\RM{o}_{\cT{Y}}\}}(\varepsilon))
    = \frac{1}{2} \sum\limits_{n=1}^{N}\abs{2c_2\zeta_f - \varepsilon + f(\RM{x}_1, \RM{v}_{q_{\delta}})^2 - u_{\cT{Y}}(\RM{v}_{q_{\delta}}, \RM{o}_{\cT{Y}})^2},
    & \hspace*{-5cm}\text{ if } q_{\delta} = M+2.
\end{cases}
\end{equation*}
Then using definition of $\BT{a}_{\cT{X}}(\varepsilon)_{\bullet, \bullet}$ (Eq.~\ref{eq:A_X}) and accounting for equation with Eq.~\ref{eq:gamma_delta} we have
\begin{equation*}
    \begin{cases}
    \abs{\gamma - \BT{c}^{\RM{x}_1}_{\cT{X} \times \{\RM{o}_{\cT{X}}\}}(\varepsilon)_{q_{\gamma}, q_{\gamma}}} \,\leq\,
    
    \frac{M-1}{2} \abs{c_1\zeta_f + c_2\zeta_f - \varepsilon}
    +
    \frac{1}{2} \abs{c_1\zeta_f + c_2\zeta_f - \varepsilon - u_{\cT{X}}(\RM{x}_1, \RM{o}_{\cT{X}})^2},
    & \hspace*{-2cm}\text{ if } \RM{v}_{q_{\gamma}} = \RM{x}_1 \text{ i.e., } q_{\gamma} = 2,
    \\
    \abs{\gamma - \BT{c}^{\RM{x}_1}_{\cT{X} \times \{\RM{o}_{\cT{X}}\}}(\varepsilon)_{q_{\gamma}, q_{\gamma}}} \,\leq\,
    \frac{M-1}{2} \abs{c_1\zeta_f + c_2\zeta_f - \varepsilon}
    +
    \frac{1}{2}  \abs{2c_2\zeta_f - \varepsilon + d_{\cT{X}}(\RM{v}_{q_{\gamma}}, \RM{x}_1)^2 - u_{\cT{X}}(\RM{x}_1, \RM{o}_{\cT{X}})^2},
    & \hspace*{-1.5cm}\text{ if } q_{\gamma} = 3, \ldots, \, M+1,
    \\
    \abs{\gamma - \BT{c}^{\RM{x}_1}_{\cT{X} \times \{\RM{o}_{\cT{X}}\}}(\varepsilon)_{q_{\gamma}, q_{\gamma}}} \,\leq\,
     
     \frac{1}{2}  \sum\limits_{m=1}^{M}\abs{2c_2\zeta_f - \varepsilon + d_{\cT{X}}(\RM{v}_{q_{\gamma}}, \RM{x}_1)^2 - u_{\cT{X}}(\RM{x}_1, \RM{o}_{\cT{X}})^2},
    & \hspace*{-2cm}\text{ if } q_{\gamma} = M+2,
    \\
    \\
    \abs{\delta - \BT{c}^{\RM{x}_1}_{\cT{Y} \times \{\RM{o}_{\cT{Y}}\}}(\varepsilon)_{q_{\delta}, q_{\delta}}} \,\leq\,
    \sum\limits_{\substack{j = M+3\\j\neq q_{\gamma}}}^{Q}
         \abs{\BT{a}^{\RM{x}_1}_{\cT{Y}}(\varepsilon)_{q_{\gamma} - (M+2), j - (M+2)}}
    +
    \frac{1}{2} \abs{2c_2\zeta_f - \varepsilon + f(\RM{x}_1, \RM{v}_{q_{\delta}})^2- u_{\cT{Y}}(\RM{v}_{q_{\delta}}, \RM{o}_{\cT{Y}})^2},
    \\
    & \hspace*{-5cm} \text{ if } q_{\delta} = M+3, \ldots,\,Q = M+N+2,
    \\

    \abs{\delta - \BT{c}^{\RM{x}_1}_{\cT{Y} \times \{\RM{o}_{\cT{Y}}\}}(\varepsilon)_{q_{\delta}, q_{\delta}}} \,\leq\,
    
    \frac{1}{2} \sum\limits_{n=1}^{N}\abs{2c_2\zeta_f - \varepsilon + f(\RM{x}_1, \RM{v}_{q_{\delta}})^2 - u_{\cT{Y}}(\RM{v}_{q_{\delta}}, \RM{o}_{\cT{Y}})^2},
    & \hspace*{-5cm}\text{ if } q_{\delta} = M+2.
\end{cases}
\end{equation*}
Using the triangle inequality, and the change of expression $c_1\zeta_f + c_2\zeta_f - \varepsilon =  2\zeta_f + c_1\zeta_f + c_2\zeta_f - \varepsilon - 2\zeta_f $ and $2c_2\zeta_f - \varepsilon = 2\zeta_f + c_1\zeta_f + c_2\zeta_f - \varepsilon + c_2\zeta_f - c_1\zeta_f - 2\zeta_f$ we have
\begin{equation*}
    \begin{cases}
    \abs{\gamma - \BT{c}^{\RM{x}_1}_{\cT{X} \times \{\RM{o}_{\cT{X}}\}}(\varepsilon)_{q_{\gamma}, q_{\gamma}}} \,\leq\,
    
    \frac{M}{2} \abs{2\zeta_f + c_1\zeta_f + c_2\zeta_f - \varepsilon}
    + M\zeta_f +
    \frac{1}{2} u_{\cT{X}}(\RM{x}_1, \RM{o}_{\cT{X}})^2,
    & \hspace*{-5cm}\text{ if } \RM{v}_{q_{\gamma}} = \RM{x}_1 \text{ i.e., } q_{\gamma} = 2,
    \\
    \abs{\gamma - \BT{c}^{\RM{x}_1}_{\cT{X} \times \{\RM{o}_{\cT{X}}\}}(\varepsilon)_{q_{\gamma}, q_{\gamma}}} \,\leq\,
    \frac{M}{2} \abs{2\zeta_f + c_1\zeta_f + c_2\zeta_f - \varepsilon}
    + M \zeta_f +
    \frac{1}{2}  \abs{c_2\zeta_f - c_1\zeta_f + d_{\cT{X}}(\RM{v}_{q_{\gamma}}, \RM{x}_1)^2 - u_{\cT{X}}(\RM{x}_1, \RM{o}_{\cT{X}})^2},
    \\
    & \hspace*{-7cm}\text{ if } q_{\gamma} = 3, \ldots, \, M+1,
    \\
    
    \abs{\gamma - \BT{c}^{\RM{x}_1}_{\cT{X} \times \{\RM{o}_{\cT{X}}\}}(\varepsilon)_{q_{\gamma}, q_{\gamma}}} \,\leq\,
     
     \frac{M}{2} \abs{2\zeta_f + c_1\zeta_f + c_2\zeta_f - \varepsilon} + M\zeta_f + \sum\limits_{m=1}^{M}\abs{c_2\zeta_f - c_1\zeta_f + d_{\cT{X}}(\RM{v}_{q_{\gamma}}, \RM{x}_1)^2 - u_{\cT{X}}(\RM{x}_1, \RM{o}_{\cT{X}})^2},
     \\
    & \hspace*{-7cm}\text{ if } q_{\gamma} = M+2,
    
    \\
    \abs{\delta - \BT{c}^{\RM{x}_1}_{\cT{Y} \times \{\RM{o}_{\cT{Y}}\}}(\varepsilon)_{q_{\delta}, q_{\delta}}} \,\leq\,
    \sum\limits_{\substack{j = M+3\\j\neq q_{\gamma}}}^{Q}
         \abs{\BT{a}^{\RM{x}_1}_{\cT{Y}}(\varepsilon)_{q_{\gamma} - (M+2), j - (M+2)}}
    +
    \frac{1}{2} \abs{2\zeta_f + c_1\zeta_f + c_2\zeta_f - \varepsilon} + \zeta_f
    \\
    & \hspace{-9cm} + \frac{1}{2} \abs{c_2\zeta_f - c_1\zeta_f + f(\RM{x}_1, \RM{v}_{q_{\delta}})^2- u_{\cT{Y}}(\RM{v}_{q_{\delta}}, \RM{o}_{\cT{Y}})^2},
    \\
    & \hspace*{-7cm} \text{ if } q_{\delta} = M+3, \ldots,\,Q = M+N+2,
    \\
  
    \abs{\delta - \BT{c}^{\RM{x}_1}_{\cT{Y} \times \{\RM{o}_{\cT{Y}}\}}(\varepsilon)_{q_{\delta}, q_{\delta}}} \,\leq\,
    
    \frac{N}{2} \abs{2\zeta_f + c_1\zeta_f + c_2\zeta_f - \varepsilon} + N\zeta_f + \sum\limits_{n=1}^{N}\abs{c_2\zeta_f - c_1\zeta_f + f(\RM{x}_1, \RM{v}_{q_{\delta}})^2 - u_{\cT{Y}}(\RM{v}_{q_{\delta}}, \RM{o}_{\cT{Y}})^2},
    & \hspace*{-0.75cm}\text{ if } q_{\delta} = M+2.
\end{cases}
\end{equation*}
Now, $\BT{c}^{\RM{x}_1}_{\cT{X} \times \{\RM{o}_{\cT{X}}\}}(\varepsilon)_{q_{\gamma}, q_{\gamma}\,/\,q_{\delta}, q_{\delta}}$ and $\BT{c}^{\RM{x}_1}_{\cT{Y} \times \{\RM{o}_{\cT{Y}}\}}(\varepsilon)_{q_{\gamma}, q_{\gamma}\,/\,q_{\delta}, q_{\delta}}$ are replaced by their respective values to get
\begin{equation*}
    \begin{cases}
    \gamma \,\geq\, \frac{c_1\zeta_f}{2}
    -
    \frac{M}{2} \abs{2\zeta_f + c_1\zeta_f + c_2\zeta_f - \varepsilon}
    - M\zeta_f -
    \frac{1}{2} u_{\cT{X}}(\RM{x}_1, \RM{o}_{\cT{X}})^2,
    & \hspace*{-5cm}\text{ if } \RM{v}_{q_{\gamma}} = \RM{x}_1 \text{ i.e., } q_{\gamma} = 2,
    \\
     \gamma \,\geq\, \frac{c_2\zeta_f}{2}
    -\frac{M}{2} \abs{2\zeta_f + c_1\zeta_f + c_2\zeta_f - \varepsilon}
    - M \zeta_f -
    \frac{1}{2}  \abs{c_2\zeta_f - c_1\zeta_f + d_{\cT{X}}(\RM{v}_{q_{\gamma}}, \RM{x}_1)^2 - u_{\cT{X}}(\RM{x}_1, \RM{o}_{\cT{X}})^2},
    \\
    & \hspace*{-7cm}\text{ if } q_{\gamma} = 3, \ldots, \, M+1,
    \\
    
    \gamma \,\geq\, \frac{c_2\zeta}{2} + f(\RM{x}_1, \RM{y}_n)^2
     -
     \frac{M}{2} \abs{2\zeta_f + c_1\zeta_f + c_2\zeta_f - \varepsilon} - M\zeta_f - \sum\limits_{m=1}^{M}\abs{c_2\zeta_f - c_1\zeta_f + d_{\cT{X}}(\RM{v}_{q_{\gamma}}, \RM{x}_1)^2 - u_{\cT{X}}(\RM{x}_1, \RM{o}_{\cT{X}})^2},
     \\
    & \hspace*{-7cm}\text{ if } q_{\gamma} = M+2,
    
    \\
    \delta \,\geq\, \frac{c_2\zeta_f}{2}
    -
    \sum\limits_{\substack{j = M+3\\j\neq q_{\gamma}}}^{Q}
         \abs{\BT{a}^{\RM{x}_1}_{\cT{Y}}(\varepsilon)_{q_{\gamma} - (M+2), j - (M+2)}}
    -
    \frac{1}{2} \abs{2\zeta_f + c_1\zeta_f + c_2\zeta_f - \varepsilon} - \zeta_f
    \\
    & \hspace{-9cm} + \frac{1}{2} \abs{c_2\zeta_f - c_1\zeta_f + f(\RM{x}_1, \RM{v}_{q_{\delta}})^2- u_{\cT{Y}}(\RM{v}_{q_{\delta}}, \RM{o}_{\cT{Y}})^2},
    \\
    & \hspace*{-7cm} \text{ if } q_{\delta} = M+3, \ldots,\,Q = M+N+2,
    \\
  
    \delta \,\geq\, \frac{c_2\zeta_f}{2}
    -
    \frac{N}{2} \abs{2\zeta_f + c_1\zeta_f + c_2\zeta_f - \varepsilon} - N\zeta_f - \sum\limits_{n=1}^{N}\abs{c_2\zeta_f - c_1\zeta_f + f(\RM{x}_1, \RM{v}_{q_{\delta}})^2 - u_{\cT{Y}}(\RM{v}_{q_{\delta}}, \RM{o}_{\cT{Y}})^2},
    & \hspace*{-0.75cm}\text{ if } q_{\delta} = M+2.
\end{cases}
\end{equation*}
Now, we replace $\BT{a}^{\RM{x}_1}_{\cT{Y}}(\varepsilon)_{\bullet, \bullet}$ by its expression and re-arranging similarly as in the previous step, we have
\begin{equation}
\begin{split}
     \sum\limits_{\substack{j = M+3\\j\neq q_{\gamma}}}^{Q}
         \abs{\BT{a}^{\RM{x}_1}_{\cT{Y}}(\varepsilon)_{q_{\gamma} - (M+2), j - (M+2)}} 
        \leq  &
    \\
    & \hspace*{-6cm}  \frac{1}{2}\underbrace{\sum\limits_{\substack{j = M+3\\j\neq q_{\gamma}}}^{Q} \left(\abs{2\zeta_f + c_1\zeta_f + c_2\zeta_f - \varepsilon} + \abs{c_2\zeta_f - c_1\zeta_f + f(\RM{x_1}, \RM{y}_{q_{\gamma} - (M+2)})^2 + f(\RM{x}_1, \RM{y}_{j - (M+2)})^2 - d_{\cT{Y}}(\RM{y}_{q_{\gamma} - (M+2)}, \RM{y}_{j - (M+2)})^2} \right)}_{\left([M-1]\abs{2\zeta_f + c_1\zeta_f + c_2\zeta_f - \varepsilon} + \sum\limits_{\substack{j = M+3\\j\neq q_{\gamma}}}^{Q} \abs{c_2\zeta_f - c_1\zeta_f + f(\RM{x_1}, \RM{y}_{q_{\gamma} - (M+2)})^2 + f(\RM{x}_1, \RM{y}_{j - (M+2)})^2 - d_{\cT{Y}}(\RM{y}_{q_{\gamma} - (M+2)}, \RM{y}_{j - (M+2)})^2}\right)}.
\end{split}
\end{equation}

Then, if we set $\varepsilon = c_3\zeta_f$, for $\gamma$ and $\delta$, we also obtain a lower bound of the form
\begin{equation*}
    \begin{cases}
    \gamma \geq c_1 - M \abs{2 + c_1+ c_2 - c_3} - K_2(c_1, c_2)\\
    \delta \geq c_2  - \max(M, N) \abs{2 + c_1 + c_2 - c_3} - \max\{K_3(c_1, c_2), K_4(c_1, c_2), K_5(c_1, c_2), K_6(c_1, c_2)\},\\
  \end{cases}
\end{equation*}
for some non-negative constants $K_i(c_1, c_2), \, i= 2,\ldots,6$ dependent on $c_1$ and $c_2$. Thus the conditions for which $\BT{C}^{\RM{x}_1}_{\cT{X} \times \{\RM{o}_{\cT{X}}\}}(\varepsilon)$ and $\BT{C}^{\RM{x}_1}_{\cT{Y} \times \{\RM{o}_{\cT{Y}}\}}(\varepsilon)$ are positive semi-definite are also equivalent to condition (i) of Theorem~\ref{thm:vareps}.
\end{proof}